\documentclass[fleqn,final]{zzz}
\usepackage{amssymb}
\usepackage{amsmath, amsfonts}
\usepackage{mathrsfs}
\usepackage{color}
\usepackage{tikz}
\usepackage{float}
\usepackage{graphicx}
\usepackage{rotating}
\usepackage[pdftex]{hyperref}
\usepackage{enumerate}
\usepackage{soul}
\usepackage{comment}
\usepackage{xparse}
\usepackage{chngcntr}
\usepackage{apptools}
\usepackage{caption}
\usepackage[normalem]{ulem}
\usepackage{subcaption}
\usepackage{arydshln}
\theoremstyle{theorem}

\theoremstyle{definition}

\numberwithin{proposition}{section}
\numberwithin{remark}{section}
\numberwithin{equation}{section}

\hypersetup{
colorlinks=true,
urlcolor=blue,
linkcolor=red,
citecolor=green
}

\NewDocumentCommand{\qfrac}{smm}{%
 \dfrac{\IfBooleanT{#1}{\vphantom{\big|}}#2}{\mathstrut #3}%
}

\makeatletter
\def\eqnarray{\stepcounter{equation}\let\@currentlabel=\theequation
\global\@eqnswtrue
\tabskip\@centering\let\\=\@eqncr
$$\halign to \displaywidth\bgroup\hfil\global\@eqcnt\z@
  $\displaystyle\tabskip\z@{##}$&\global\@eqcnt\@ne
  \hfil$\displaystyle{{}##{}}$\hfil
  &\global\@eqcnt\tw@ $\displaystyle{##}$\hfil
  \tabskip\@centering&\llap{##}\tabskip\z@\cr}

\def\endeqnarray{\@@eqncr\egroup
      \global\advance\c@equation\m@ne$$\global\@ignoretrue}

\def\@yeqncr{\@ifnextchar [{\@xeqncr}{\@xeqncr[5pt]}}
\makeatother

\newcommand{\R}{{\mathbb R}}

\definecolor{darkgreen}{rgb}{0.0, 0.21, 0.06}

\begin{document}

\renewcommand{\PaperNumber}{***}

\FirstPageHeading

\ShortArticleName{%
%Symmetry of terminating basic hypergeometric representations of the Askey--Wilson Polynomials
Closed form for $\sum_{k=1}^n k^p$ through the Hermite integral representation...
}

\ArticleName{%
%Symmetry of terminating basic hypergeometric\\ representations of the Askey--Wilson polynomials
Closed form for $\sum_{k=1}^n k^p$ through the Hermite integral representation of the Hurwitz zeta function
}

% Names of the authors for the title of the paper
\Author{Abdulhafeez A. Abdulsalam\,$^\dag\!\!\ $}

\AuthorNameForHeading{A.~A.~Abdulsalam}
\Address{$^\dag$ Department of Mathematics, University of Ibadan, Ibadan, Oyo, Nigeria}
\EmailD{hafeez147258369@gmail.com} % E-mail address of First Author
% E-mail address of First Author

%\Address{$^{\S\S}$ Department of Mathematics,
%University of Rochester, Rochester, NY 14627, USA
%% Address of First Author
%}
%\EmailD{random9483@gmail.com} % E-mail address of First Author

\ArticleDates{Received \today~in final form ????; 
Published online ????}

\Abstract{Sums of positive integer powers have captivated the attention of mathematicians since ancient times. Over the centuries, mathematicians from diverse backgrounds have provided expressions for the sum of positive integer powers of the first $n$ positive integers. In this paper, we contribute to this endeavour by deriving Bernoulli's formula for  $\sum_{k=1}^n k^p,$ for all positive integers $n$ and nonegative integers $p$, through the utilization of the Hermite integral representation of the Hurwitz zeta function.
}
%.}

%Please type here List of Keywords for your article separated 
%by semicolon.
% Keywords required only for MST, PB, PMB, PM, JOA, JOB?
% Keywords:

\Classification{11B68, 11M06, 11M35, 33B15}
%{??????} % e.g. 35A30; 81Q05
%For 2010 Mathematics Subject Classification see
%http://www.ams.org/mathscinet/msc/msc2010.html
\allowdisplaybreaks
\section{Introduction}
The sum $S_k(n)$, defined as
\[
S_p(n) = \sum_{k=1}^n k^p,
\]
for all positive integers $n$ and nonegative integers $p$, represents the sum of the $p$th powers of the first $n$ integers. In this paper, we prove Bernoulli's generalization for $S_p(n)$ using the Hermite integral representation of the Hurwitz zeta function. In Beery's work, historical accounts are meticulously detailed. Beery \cite{bib4} delves into the historical context of the problem, tracing the emergence of generalizable formulas for $S_p(n)$ back to the contributions of notable mathematicians. Notably, the West first saw generalizable formulas for sums of integer powers introduced by Thomas Harriot (1560-1621) of England. Around the same period, Johann Faulhaber (1580-1635) of Germany presented formulas for these sums up to the 17th power, surpassing previous efforts in terms of power level but lacking clarity in generalization. While Pierre de Fermat (1601-1665) is often credited with discovering these formulas, it was his fellow French mathematician Blaise Pascal (1623-1662) who offered more explicit formulations. However, it was the Swiss mathematician Jakob Bernoulli (1654-1705) who gained recognition for presenting the most useful and generalizable formulation within the European mathematical community. His formula for the sum \(S_p(n)\) is as follows \cite{bib4}
\begin{equation}\label{bernb1}
S_p(n) = \frac{1}{p+1} \sum_{k=0}^p \binom{p+1}{k} B_k n^{p-k+1},
\end{equation}
where $B_k$ is the $k$th Bernoulli number  \cite[(3)]{bib5} and 
\[\frac{1}{p+1} \sum_{k=0}^p \binom{p+1}{k} B_k = 1.\]
Bernoulli's work marked a significant milestone in the development of these formulas and remains a pivotal reference point in the history of this problem. We remark that formula \eqref{bernb1} assumes $B_1 = \frac{1}{2}$, in contrast to the widely used value of $B_1$ as $B_1 = -\frac{1}{2}$. The value $B_1 = -\frac{1}{2}$ is used in NIST's Digital Library of Mathematical Functions \cite[\S 24.2(iv)]{bib1}, \textsf{Mathematica 13.0} and \textsf{Maple 2022}. Noting this observation, utilizing formula \eqref{bernb1}, we derive the following first four terms of $S_p(n)$ 
\[\sum_{k=1}^n k^0 = n, \quad \sum_{k=1}^n k = \frac{n}{2}(n + 1), \quad \sum_{k=1}^n k^2 = \frac{n}{6}(2n^2 + 3n +1), \quad \sum_{k=1}^n k^3 = \frac{n^2(n+1)^2}{4}.\]
The Bernoulli polynomials $B_n(x)$ are defined by \cite[(24.2.3)]{bib1}
\begin{equation}\label{multiv1}
\frac{t e^{xt}}{e^t - 1} = \sum_{n=0}^\infty  \frac{B_n(x) t^n}{n!},
\end{equation}
where $B_n = B_n(0)$ \cite[(24.2.4)]{bib1}. The discussion regarding the value of $B_1$ stems from the definition provided in \eqref{multiv1}. For if we set $x=1$ in \eqref{multiv1}, we obtain
\begin{equation}\label{multiv2}
\frac{te^t}{e^{t} - 1} = \sum_{n=0}^\infty \frac{B_n(1) t^n}{n!}.
\end{equation}
By expanding the series on the left-hand side of \eqref{multiv2} about $t=0$, we have
\begin{equation}\label{multiv3}
1 + \frac{t}{2} + \frac{t^2}{12} - \frac{t^4}{720} + O(t^4) = \sum_{n=0}^\infty \frac{B_n(1) t^n}{n!}.
\end{equation}
Upon comparison, we find that $B_1(1) = \frac{1}{2}$. Now, by setting $x=0$ in \eqref{multiv1}, we have
\[\frac{t}{e^t - 1} = \sum_{n=0}^\infty \frac{B_n t^n}{n!}.\]
Expanding the series on the left-hand side of \eqref{multiv3} about $t=0$, we get
\[1 - \frac{t}{2} + \frac{t^2}{12} - \frac{t^4}{720} + O(t^4) = \sum_{n=0}^\infty \frac{B_n t^n}{n!}.\]
Upon comparison, we find that $B_1 = -\frac{1}{2}$. Notice that the two series differ only in the coefficient of $ t $, as $\frac{te^t}{e^t - 1} = t + \frac{t}{e^t - 1}$. In particular, this distinction is recognized in the literature, where $B_1^- := B_1 = -\frac{1}{2}$ and $B_1^+ = B_1(1) = \frac{1}{2}$, and $B_n(1) = B_n$ for all nonnegative integers $n \neq 1$. When calculating the area under the curve $y = x^p$ over the interval $[0, 1]$ using the Riemann sum, one has
\[
\int_0^1 x^p \, \mathrm{d}x = \lim_{n \to \infty} \sum_{k=1}^n \frac{1}{n} \left(\frac{k}{n}\right)^p = \lim_{n \to \infty} \frac{1}{n^{p+1}} \left(\sum_{k=1}^n k^p \right).
\]
It is evident that we need not consider the sum $\sum_{k=1}^n k^p$ to conclude that the area of the under the curve $y = x^p$ over $[0, 1]$ is always $\frac{1}{p+1}$. However, if one wishes to calculate the area under the curve using the Riemann sum, it is essential to have an explicit expression for the sum $S_p(n)$ for any positive integer $p$. The simplest known recursive formula for the generalized sum $S_p(n)$ is \cite[(11)]{bib2}
\[
(1 + n)^{k+1} = 1 + \sum_{r=0}^p \binom{k+1}{r} S_r(n),
\]
where integers $p \geq 0$, $n > 0$, and $S_p(n)$ could be deduced from $S_0(n), S_1(n), \ldots, S_{p-1}(n)$. Additional recursive formulas are provided in \cite{bib2}. Schultz \cite{bib3} noted that standard closed forms for $S_p(n)$ typically involve Bernoulli numbers or Stirling numbers of the second kind, which are derived from complex recurrence relations. One example is Bernoulli's formula \eqref{bernb1}, which features the Bernoulli numbers. Another formula for $S_p(n)$, not explicitly stated by Bernoulli, involving the Bernoulli numbers is given by \cite[(1)]{bib5}
\[S_p(n) = \frac{B_{p+1}(n) - B_{p+1}}{p+1},\]
where integers $n \geq 2$, $p \geq 1$, and $B_n(x)$ the Bernoulli polynomial \cite[(2)]{bib5}. In an effort to make the derivation of a formula for $S_p(n)$ less complicated, we derive a formula for $S_p(n)$ using the Hermite integral representation of the Hurwitz zeta function. The Riemann zeta function \cite[\S 25.2]{bib1} and the Hurwitz zeta function \cite[\S 25.11]{bib1} are, respectively, defined as
$$\zeta(s) = \sum_{n=1}^{\infty} \frac{1}{n^s}, \quad \zeta(s, z) = \sum_{n=0}^{\infty} \frac{1}{(n+z)^s},$$
where $z \not\in -\mathbb{N}_0$, $\Re\,s > 1$. The domain $\Re\,s > 1$ can be extended to $s \in \mathbb{C}\setminus \{1\}$ through analytic continuation, using for instance, the Hermite integral representation for the Hurwitz zeta function \cite[(25.11.29)]{bib1}
\begin{equation}\label{herm}
\zeta(s, z) = \frac{z^{-s}}{2} + \frac{z^{1-s}}{s-1} + 2\int_0^\infty \frac{\sin\left(s\arctan\left(x/z\right)\right)}{\left(x^2 + z^2\right)^{\frac{s}{2}}\left(\mathrm{e}^{2\pi x} - 1\right)}\, \mathrm{d}x.
\end{equation}
Throughout this article, $\Im z$ denotes the imaginary part of $z$.
\section{Proof using the integral representation of the Hurwitz zeta function}
In this section, we establish Bernoulli's formula \eqref{bernb1} for $S_p(n)$ using Hermite's integral representation of the Hurwitz zeta function. Two theorems concerning $S_p(n)$ are presented, arising from the value assumed for $B_1$.
\begin{theorem} Suppose $B_1 = -\frac{1}{2}$ and let $n$ and $p$ be integers such that $n \geq 1$, $p \geq 0$. Then \label{thm1} 
\begin{equation}\label{mainthh1}
\begin{split}
S_p(n) &= n, \quad \textup{if $p = 0$},
\\&= \begin{cases}  \displaystyle \frac{1}{p+1} \sum_{k=0}^{\frac{p}{2}} \binom{p+1}{2k+1} (n+1)^{2k+1}B_{p-2k} - \frac{(n+1)^{p}}{2}, &\textup{if $p$ is even,}  \\ \displaystyle \frac{1}{p+1} \sum_{k=1}^{\frac{p+1}{2}} \binom{p+1}{2k} (n+1)^{2k} B_{p-2k+1} - \frac{(n+1)^{p}}{2},  & \textup{if $p$ is odd.} \end{cases}
\\&= \frac{1}{p+1} \sum_{k=0}^{p} \binom{p+1}{k} n^{p-k+1} B_{k} + n^p = \frac{1}{p+1} \sum_{k=0}^{p} (-1)^k \binom{p+1}{k} n^{p-k+1} B_{k},
\end{split}
\end{equation}
where $B_k$ is the $k$th Bernoulli number.
\end{theorem}

\begin{proof}
Expressing $S_p(n)$ in terms of the zeta functions, we have
\begin{equation}\label{splhaf1}
\begin{split}
S_p(n) &= \sum_{k=1}^\infty k^p - \sum_{k=n+1}^\infty k^p
\\&= \sum_{k=1}^\infty k^p - \sum_{k=0}^\infty (k+n+1)^p = \zeta(-p) - \zeta(-p, n+1).
\end{split}
\end{equation}
Since the argument $p$ in $\zeta(-p,n+1)$ is positive, we use the Hermite integral representation, as it provides analytic continuation of the domain $\Re\,s > 1$ of the Hurwitz zeta function $\zeta(s, z)$ to $s \in \mathbb{C}\setminus \{1\}$. This yields
\begin{equation}\label{lastinth1}
\zeta(-p, n+1) = \frac{(n+1)^{p}}{2} - \frac{(n+1)^{1+p}}{1+p} - 2\int_0^\infty \frac{\sin\left(p\arctan\left(x/(n+1)\right)\right)}{\left(x^2 + (n+1)^2\right)^{\frac{-p}{2}}\left(\mathrm{e}^{2\pi x} - 1\right)}\, \mathrm{d}x.
\end{equation}
Define the last integral in \eqref{lastinth1} as
\[f(p, n) = \int_0^\infty \frac{\sin\left(p\arctan\left(x/(n+1)\right)\right)}{\left(x^2 + (n+1)^2\right)^{\frac{-p}{2}}\left(\mathrm{e}^{2\pi x} - 1\right)}\, \mathrm{d}x.\]
By De Moivre's theorem \cite[(4.21.34)]{bib1}, we have
\[f(p, n) = \Im \int_0^\infty \frac{\left(\cos\left(\arctan\left(\frac{x}{n+1}\right)\right) + \mathrm{i}\sin\left(\arctan\left(\frac{x}{n+1}\right)\right)\right)^p}{\left(x^2 + (n+1)^2\right)^{\frac{-p}{2}}\left(\mathrm{e}^{2\pi x} - 1\right)}\, \mathrm{d}x,\]
where $\mathrm{i}=\sqrt{-1}$. By the binomial theorem, we have
\[f(p, n) = \Im \int_0^\infty\sum_{m=0}^p \binom{p}{m} \frac{\left(\cos\left(\arctan\left(\frac{x}{n+1}\right)\right)\right)^m \mathrm{i}^{p-m} \left(\sin\left(\arctan\left(\frac{x}{n+1}\right)\right)\right)^{p-m}}{\left(x^2 + (n+1)^2\right)^{\frac{-p}{2}}\left(\mathrm{e}^{2\pi x} - 1\right)}\, \mathrm{d}x.\]
Interchanging summation and integration, we have
\[f(p, n) = \Im \sum_{m=0}^p \binom{p}{m} \int_0^\infty\frac{\left(\cos\left(\arctan\left(\frac{x}{n+1}\right)\right)\right)^m \mathrm{i}^{p-m} \left(\sin\left(\arctan\left(\frac{x}{n+1}\right)\right)\right)^{p-m}}{\left(x^2 + (n+1)^2\right)^{\frac{-p}{2}}\left(\mathrm{e}^{2\pi x} - 1\right)}\, \mathrm{d}x.\]
Using the trigonometric identity $1 + \tan^2\theta = \sec^2\theta$, $\theta \in \R$, we deduce
\[\cos\left(\arctan\left(\frac{x}{n+1}\right)\right) = \frac{1}{\sqrt{1 + \tan^2\left(\arctan\left(\frac{x}{n+1}\right)\right)}} = \frac{n+1}{\sqrt{x^2 + (n+1)^2}},\]
\[\sin\left(\arctan\left(\frac{x}{n+1}\right)\right) = \sqrt{1-\frac{(n+1)^2}{x^2 + (n+1)^2}} = \frac{x}{\sqrt{x^2 + (n+1)^2}}.\]
Therefore, we can express $f(p, n)$ as
\begin{align*}
f(p, n) &= \Im \sum_{m=0}^p \binom{p}{m} (n+1)^m \mathrm{i}^{p-m}\int_0^\infty \frac{x^{p-m}}{\left(x^2 + (n+1)^2\right)^{\frac{p}{2}} \left(x^2 + (n+1)^2\right)^{\frac{-p}{2}}\left(\mathrm{e}^{2\pi x} - 1\right)}\, \mathrm{d}x
\\&= \Im \sum_{m=0}^p \binom{p}{m} (n+1)^m \mathrm{i}^{p-m}\int_0^\infty \frac{x^{p-m}}{\mathrm{e}^{2\pi x} - 1}\, \mathrm{d}x 
\\&= \Im \sum_{m=0}^p \binom{p}{m} (n+1)^m \mathrm{i}^{p-m} \sum_{l=0}^\infty \int_0^\infty  x^{p-m} \mathrm{e}^{-2\pi x (l+1)} \, \mathrm{d}x 
\\&= \Im \sum_{m=0}^p \binom{p}{m} (n+1)^m \mathrm{i}^{p-m} \sum_{l=0}^\infty \frac{(p-m)!}{(2\pi)^{p-m+1} (l+1)^{p-m+1}} 
\\&= \Im \sum_{m=0}^p \binom{p}{m} (n+1)^m \mathrm{i}^{p-m}\frac{(p-m)!}{(2\pi)^{p-m+1}} \zeta(p-m+1).
\end{align*}
Since
\[\Im  \mathrm{i}^{p-m} = \Im e^{\mathrm{i}\frac{\pi}{2}(p-m)} = \Im\left(\cos\left(\frac{\pi}{2}(p-m)\right) +\mathrm{i}\sin\left(\frac{\pi}{2}(p-m)\right)\right) = \sin\left(\frac{\pi}{2}(p-m)\right),\]
we conclude that
\begin{equation}\label{sether1}
f(p, n) = \sum_{m=0}^p  \frac{p! (n+1)^m}{m! (2\pi)^{p-m+1}} \sin\left(\frac{\pi}{2}(p-m)\right)\zeta(p-m+1).
\end{equation}
Taking into consideration that $\zeta(-p) = -\frac{B_{p+1}}{p+1}$ for $p \neq 0$ \cite[(25.11.14)]{bib1}, and the fact that $\sin\left(\frac{\pi}{2}(p-p)\right) = 0$, we can substitute \eqref{sether1} into \eqref{lastinth1} and subsequently substitute the latter into \eqref{splhaf1}. This allows us to express $S_p(n)$ as
\begin{equation}\label{pkts}
\begin{split}
S_p(n) &= -\frac{B_{p +1}}{p+1} - \frac{(n+1)^{p}}{2} + \frac{(n+1)^{1+p}}{1+p} 
\\&\quad+ 2\sum_{m=0}^{p-1} \frac{p! (n+1)^m}{m! (2\pi)^{p-m+1}}\sin\left(\frac{\pi}{2}(p-m)\right)\zeta(p-m+1), 
\end{split}
\end{equation}
where $p$ is a positive integer. The last sum in \eqref{pkts} can be further simplified for even and odd cases of $p$. To this end, by setting $p = 2\rho-1$ in \eqref{sether1}, where $\rho$ is a positive integer, we have
\[f(2\rho-1, n) = \sum_{m=0}^{2\rho-2} \frac{(2\rho-1)! (n+1)^m}{m! (2\pi)^{2\rho-m}}\sin\left(\frac{\pi}{2}(2\rho-m-1)\right)\zeta(2\rho-m).\]
When $m$ is odd, $\sin\left(\frac{\pi}{2}(2\rho-m-1)\right) = 0$. This implies
\begin{align*}
f(2\rho-1, n) &= \sum_{m=0}^{\rho-1} \frac{(2\rho-1)! (n+1)^{2m}}{(2m)! (2\pi)^{2\rho-2m}}\sin\left(\frac{\pi}{2}(2\rho-2m-1)\right)\zeta(2\rho-2m).
\end{align*}
As $\sin\left(\frac{\pi}{2}(2\rho-2m-1)\right) = (-1)^{p-m+1}$ and \cite[(7)]{bib5}
\[\zeta(2\rho-2m) = (-1)^{p-m-1} \frac{ (2\pi)^{2\rho-2m} B_{2\rho-2m}}{2(2\rho-2m)!},\]
we derive
\begin{equation}\label{efoddrh}
f(2\rho-1, n) = \frac{(2\rho-1)!}{2}\sum_{m=0}^{\rho-1} \frac{(n+1)^{2m}}{(2m)! (2\rho-2m)!} B_{2\rho-2m}.
\end{equation}
Utilizing \eqref{efoddrh} in \eqref{pkts}, we deduce
\begin{equation}\label{threerho1}
\begin{split}
S_{2\rho-1}(n) &= -\frac{B_{2\rho}}{2\rho} - \frac{(n+1)^{2\rho-1}}{2} + \frac{(n+1)^{2\rho}}{2\rho} + (2\rho-1)! \sum_{m=0}^{\rho-1} \frac{(n+1)^{2m}}{(2m)! (2\rho-2m)!} B_{2\rho-2m}
\\&=  (2\rho-1)! \sum_{m=0}^{\rho} \frac{(n+1)^{2m}}{(2m)! (2\rho-2m)!} B_{2\rho-2m} -\frac{B_{2\rho}}{2\rho} - \frac{(n+1)^{2\rho-1}}{2} 
\\&= (2\rho-1)! \sum_{m=1}^{\rho} \frac{(n+1)^{2m}}{(2m)! (2\rho-2m)!} B_{2\rho-2m}  - \frac{(n+1)^{2\rho-1}}{2}.
\end{split}
\end{equation}
By setting $k = 2\rho$ in \eqref{sether1}, where $\rho$ is a positive integer, we have
\[f(2\rho, n) = \sum_{m=0}^{2\rho-1} \frac{(2\rho)! (n+1)^m}{m! (2\pi)^{2\rho-m+1}}\sin\left(\frac{\pi}{2}(2\rho-m)\right)\zeta(2\rho-m+1).\]
When $m$ is even, $\sin\left(\frac{\pi}{2}(2\rho-m)\right) = 0$. This implies
\[f(2\rho, n) = \sum_{m=0}^{\rho-1} \frac{(2\rho)! (n+1)^{2m+1}}{(2m+1)! (2\pi)^{2\rho-2m}}\sin\left(\frac{\pi}{2}(2\rho-2m-1)\right)\zeta(2\rho-2m),\]
and so,
\begin{equation}\label{thray}
f(2\rho, n) = \frac{(2\rho)!}{2}\sum_{m=0}^{\rho-1} \frac{(n+1)^{2m+1}}{(2m+1)! (2\rho-2m)!} B_{2\rho-2m}.
\end{equation}
Utilizing \eqref{thray} in \eqref{pkts}, we deduce
\[S_{2\rho}(n) = -\frac{B_{2\rho+1}}{2\rho+1} - \frac{(n+1)^{2\rho}}{2} + \frac{(n+1)^{2\rho+1}}{2\rho+1} + (2\rho)! \sum_{m=0}^{\rho-1} \frac{(n+1)^{2m_1}  B_{2\rho-2m}}{(2m+1)! (2\rho-2m)!}.\]
Since $B_{2\rho+1} = 0$ \cite[p.~179]{bib5} and $B_0 = 1$ \cite[(3)]{bib5}, we have
\begin{equation}\label{tworho1}
\begin{split}
S_{2\rho}(n) &=  (2\rho)! \sum_{m=0}^{\rho-1} \frac{(n+1)^{2m+1} B_{2\rho-2m} }{(2m+1)! (2\rho-2m)!} + \frac{(n+1)^{2\rho+1}}{2\rho+1} - \frac{(n+1)^{2\rho}}{2}
\\&= (2\rho)! \sum_{m=0}^{\rho} \frac{(n+1)^{2m+1}B_{2\rho-2m}}{(2m+1)! (2\rho-2m)!}  - \frac{(n+1)^{2\rho}}{2}.
\end{split}
\end{equation}
By Setting $2\rho = p$ in \eqref{threerho1} and \eqref{tworho1}, we can express the formulas for $S_k(n)$ as
\begin{equation}\label{shrtthm1}
\sum_{k=1}^n k^p = \begin{cases}  \displaystyle p! \sum_{k=0}^{\frac{p}{2}} \frac{(n+1)^{2k+1}B_{p-2k}}{(2k+1)! (p-2k)!}  - \frac{(n+1)^{p}}{2}, & \textup{if $p$ is even}, \\ \\ \displaystyle p! \sum_{k=1}^{\frac{p+1}{2}} \frac{(n+1)^{2k} B_{p-2k+1}}{(2k)! (p-2k+1)!} - \frac{(n+1)^{p}}{2},  & \textup{if $p$ is odd.} \end{cases}
\end{equation}
Employing the combination notation, we can restate \eqref{shrtthm1} as
\begin{equation}\label{thislast}
S_p(n) = \sum_{k=1}^n k^p = \begin{cases}  \displaystyle \frac{1}{p+1} \sum_{k=0}^{\frac{p}{2}} \binom{p+1}{2k+1} (n+1)^{2k+1}B_{p-2k} - \frac{(n+1)^{p}}{2}, & \textup{if $p$ is even}, \\ \\ \displaystyle \frac{1}{p+1} \sum_{k=1}^{\frac{p+1}{2}} \binom{p+1}{2k} (n+1)^{2k} B_{p-2k+1} - \frac{(n+1)^{p}}{2},  & \textup{if $p$ is odd.} \end{cases}
\end{equation}
We notice the similarity between equations \eqref{bernb1} and \eqref{thislast}. To explore the possibility that they coincide, we express equation \eqref{mainthh1} in the following manner. By substituting $n$ with $n-1$ in \eqref{thislast}, we get
\[S_p(n-1)= \begin{cases}  \displaystyle \frac{1}{p+1} \sum_{k=0}^{\frac{p}{2}} \binom{p+1}{2k+1} n^{2k+1}B_{p-2k} - \frac{n^{p}}{2}, & \textup{if $p$ is even}, \\ \\ \displaystyle \frac{1}{p+1} \sum_{k=1}^{\frac{p+1}{2}} \binom{p+1}{2k} n^{2k} B_{p-2k+1} - \frac{n^{p}}{2},  & \textup{if $p$ is odd.} \end{cases}\]
Adding $n^p$ to both sides, we obtain
\[S_p(n) = \begin{cases}  \displaystyle \frac{1}{p+1} \sum_{k=0}^{\frac{p}{2}} \binom{p+1}{2k+1} n^{2k+1}B_{p-2k} + \frac{n^{p}}{2}, & \textup{if $p$ is even}, \\ \\ \displaystyle \frac{1}{p+1} \sum_{k=1}^{\frac{p+1}{2}} \binom{p+1}{2k} n^{2k} B_{p-2k+1} + \frac{n^{p}}{2},  & \textup{if $p$ is odd.} \end{cases}\]
Exploiting the well-known property that $\binom{n}{r} = \binom{n}{n-r}$, we can express $S_p(n)$ as
\begin{equation}\label{youfinally}
S_p(n) = \begin{cases}  \displaystyle \frac{1}{p+1} \sum_{k=0}^{\frac{p}{2}} \binom{p+1}{p-2k} n^{2k+1}B_{p-2k} + \frac{n^{p}}{2}, & \textup{if $p$ is even}, \\ \\ \displaystyle \frac{1}{p+1} \sum_{k=1}^{\frac{p+1}{2}} \binom{p+1}{p-2k+1} n^{2k} B_{p-2k+1} + \frac{n^{p}}{2},  & \textup{if $p$ is odd.} \end{cases}
\end{equation}
Since $B_{p-2k+1} = 0$, for $p > 2k$, $p$ even, and $B_{p-2k} = 0$, for $p > 2k+1$, $p$ odd, we can write
\begin{equation}\label{pof1}
S_p(n) = \begin{cases}  \begin{aligned}[t] &\frac{1}{p+1} \sum_{k=0}^{\frac{p}{2}} \binom{p+1}{p-2k} n^{2k+1}B_{p-2k}  \\&\quad+  \frac{1}{p+1} \sum_{k=1}^{\frac{p}{2}} \binom{p+1}{p-2k+1} n^{2k}B_{p-2k+1} + n^p,\end{aligned} & \textup{if $p$ is even}, \\ \\ \begin{aligned}[t]  & \frac{1}{p+1} \sum_{k=1}^{\frac{p+1}{2}} \binom{p+1}{p-2k+1} n^{2k} B_{p-2k+1} \\&\quad + \frac{1}{p+1} \sum_{k=0}^{\frac{p-1}{2}} \binom{p+1}{p-2k} n^{2k+1} B_{p-2k} + n^p,\end{aligned}  & \textup{if $p$ is odd.} \end{cases}
\end{equation}
If we define two functions $\alpha_p(n, k)$ and $\beta_p(n, k)$ such that
\[\alpha_p(n, k) = \frac{1}{p+1} \binom{p+1}{p-2k+1} n^{2k}B_{p-2k+1}, \quad \beta_p(n, k) = \frac{1}{p+1}  \binom{p+1}{p-2k} n^{2k+1} B_{p-2k},\]
then the reason for the last term in \eqref{pof1} being $n^p$ is due to the fact that at $k=\frac{p}{2}$, for $p$ even 
\begin{equation}\label{toberef1}
\alpha_p\left(n, \frac{p}{2}\right) = n^p B_1 = -\frac{n^p}{2},
\end{equation}
and at $k=\frac{p-1}{2}$, for $p$ odd
\begin{equation}\label{toberef2}
\beta_p\left(n, \frac{p-1}{2}\right) = n^p B_1 = -\frac{n^p}{2}.
\end{equation}
This nonzero term in the series $\sum_{k=1}^{\frac{p}{2}} \alpha_p(n, k)$ and $\sum_{k=0}^{\frac{p-1}{2}} \beta_p(n, k)$ requires us to add $\frac{n^p}{2}$, so that we can have 
\[\sum_{k=1}^{\frac{p}{2}} \alpha_p(n, k) + \frac{n^p}{2} = 0, \quad \sum_{k=0}^{\frac{p-1}{2}} \beta_p(n, k)  + \frac{n^p}{2} = 0.\]
Equation \eqref{pof1} allows us to simplify $S_p(n)$ to
\begin{equation}\label{toberef3}
S_p(n) = \frac{1}{p+1} \sum_{k=0}^{p} \binom{p+1}{p-k} n^{k+1} B_{p-k} + n^p,
\end{equation}
for all positive integers $p$. By replacing the dummy $k$ with $p - k$ in \eqref{toberef3}, we establish that for all positive integers $p$
\begin{equation}\label{toberef3hn}
S_p(n) = \frac{1}{p+1} \sum_{k=0}^{p}  \binom{p+1}{k} n^{p-k+1} B_{k} + n^p.
\end{equation}
Noting that $\zeta(0) = -\frac{1}{2}$ \cite[(25.6.1)]{bib1} in \eqref{splhaf1}, and $f(0, n) = 0$ in \eqref{sether1}, we leverage \eqref{splhaf1} and \eqref{lastinth1} to obtain
\[S_0(n) = -\frac{1}{2} - \frac{1}{2} + n+1 = n.\]
By replacing $B_n$ with $(-1)^n B_n$ in the sum in \eqref{toberef3hn}, we note that this change does not affect the sum, except at $n=1$, since $(-1)^n B_n = B_n$ for even $n$, and $(-1)^n B_n = B_n$ for odd $n \neq 1$, given that $B_{n} = 0$ for odd $n$. Thus, we derive
\[
\sum_{k=0}^{p}  \binom{p+1}{k} n^{p-k+1} B_{k}  = \sum_{k=0}^{p} (-1)^k \binom{p+1}{k} n^{p-k+1} B_{k} - n^p.
\]
This concludes the proof of the Theorem \ref{thm1} for all nonnegative integers $p$.
\end{proof}

\begin{theorem} Suppose $B_1 = \frac{1}{2}$ and let $n$ and $p$ be integers such that $n \geq 1$, $p \geq 0$. Then \label{thm2}
\begin{equation}\label{mainthh1}
\begin{split}
S_p(n) &= n, \quad \textup{if $p = 0$},
\\&= \begin{cases}  \displaystyle \frac{1}{p+1} \sum_{k=0}^{\frac{p}{2}} \binom{p+1}{2k+1} (n+1)^{2k+1}B_{p-2k} - \frac{(n+1)^{p}}{2},  &\textup{if $p$ is even,}  \\ \displaystyle \frac{1}{p+1} \sum_{k=1}^{\frac{p+1}{2}} \binom{p+1}{2k} (n+1)^{2k} B_{p-2k+1} - \frac{(n+1)^{p}}{2},  & \textup{if $p$ is odd.} \end{cases}
\\&= \frac{1}{p+1} \sum_{k=0}^{p} \binom{p+1}{k} n^{p-k+1} B_{k},
\end{split}
\end{equation}
where $B_k$ is the $k$th Bernoulli number.
\end{theorem}

\begin{proof}
The proof of Theorem \ref{thm2} follows from substituting the value $B_1 = \frac{1}{2}$ into both \eqref{toberef1} and \eqref{toberef2}, and then employing the resulting expressions in \eqref{youfinally}.
\end{proof}

\begin{remark}
If $B_1$ is assumed to be $\pm \frac{1}{2}$, then it has no effect on the formula
\[S_p(n) = \begin{cases}  \displaystyle \frac{1}{p+1} \sum_{k=0}^{\frac{p}{2}} \binom{p+1}{2k+1} (n+1)^{2k+1}B_{p-2k} - \frac{(n+1)^{p}}{2},  &\textup{if $p$ is even,}  \\ \displaystyle \frac{1}{p+1} \sum_{k=1}^{\frac{p+1}{2}} \binom{p+1}{2k} (n+1)^{2k} B_{p-2k+1} - \frac{(n+1)^{p}}{2},  & \textup{if $p$ is odd,} \end{cases}\]
as $B_{p-2k} \neq B_1$, for $p$ even and $p > 2k$, and $B_{p-2k+1} \neq B_1$, for $p$ odd and $p > 2k-1$.
\end{remark}
\section{Conclusion}
We have successfully established Bernoulli's formula for $S_p(n)$, utilizing the Hermite integral representation of the Hurwitz zeta function. Additionally, we presented two theorems regarding $S_p(n)$, depending on which sign of $\frac{1}{2}$ is assumed for the Bernoulli number $B_1$, and we note a formula unaffected by this assumption. These formulas bear particular significance due to their capacity to expedite the computation of sums involving any positive integer powers of the first $n$ positive integers. This heightened computational efficiency provides a valuable advantage when working with such summations. For instance, it is widely known that Bernoulli proudly announced that, with the assistance of the final entry in his list of formulas, he calculated the sum of the tenth powers of the first 1000 numbers in less than a quarter of an hour, resulting in the sum 91,409,924,241,424,243,424,241,924,242,500.

\section*{Acknowledgment}
I would like to express my gratitude to the Spirit of Ramanujan (SOR) STEM Talent Initiative, directed by Ken Ono, for providing me with the computational tools that were instrumental in verifying and validating my results. Special thanks are due to my undergraduate project supervisor at the University of Ibadan, Dr.~H.~P.~Adeyemo, for his invaluable introduction to the concept of power sums.
\section*{Funding}
The author did not receive funding from any organization for the submitted work.
\bibliographystyle{unsrt}
\bibliography{hermite_application}

\end{document}